\documentclass[titlepage,twoside,12pt]{article}
\usepackage{amsmath}
\usepackage{amssymb}
\usepackage{amsfonts}
\begin{document}

{\LARGE \bf Local Functions : Algebras, Ideals, \\ \\ and Reduced Power Algebras} \\ \\

{\bf Elem\'{e}r E ~Rosinger} \\ \\
{\small \it Department of Mathematics \\ and Applied Mathematics} \\
{\small \it University of Pretoria} \\
{\small \it Pretoria} \\
{\small \it 0002 South Africa} \\
{\small \it eerosinger@hotmail.com} \\ \\

\hfill {\it Dedicated to Marie-Louise Nykamp} \\ \\

{\bf Abstract} \\

A further significant extension is presented of the infinitely large class of differential algebras of generalized functions which are the basic structures in the nonlinear algebraic theory listed under 46F30 in the AMS Mathematical Subject Classification. These algebras are constructed as {\it reduced powers}, when seen in terms of Model Theory. The major advantage of these differential algebras of generalized functions is that they allow their elements to have singularities on {\it dense} subsets of their domain of definition, and {\it without} any restrictions on the respective generalized functions in the neighbourhood of their singularities. Their applications have so far been in 1) solving large classes of systems of nonlinear PDEs, 2) highly singular problems in Differential Geometry, with respective applications in modern Physics, including General Relativity and Quantum Gravity. These infinite classes of algebras contain as a particular case the Colombeau algebras, since in the latter algebras rather strongly limiting growth conditions, namely, of polynomial type, are required on the generalized functions in the neighbourhood of their singularities.  \\ \\ \\

\hspace*{3cm} ''We do not possess any method at all to derive \\
\hspace*{3cm} systematically solutions that are free of \\
\hspace*{3cm} singularities...'' \\
\hspace*{3cm} Albert Einstein : The Meaning of Relativity. \\
\hspace*{3cm} Princeton Univ. Press, 1956, p. 165 \\ \\

{\bf 0. Preliminaries} \\

The nonlinear theory of generalized functions, see 46F30 in the AMS Mathematical Subject Classification of the AMS, has known a wide range of applications in solving large classes of nonlinear PDEs, see [6-8,10-12,15-17,20,22,33,34,36,40,44,54] and the literature cited there, as well as in Abstract Differential Geometry with applications to modern Physics, including General Relativity and Quantum Gravity, [9,13-19,21,56-61], as well as in the study of manifolds, [31]. \\

These algebras are constructed as {\it reduced powers}, when seen in terms of Model Theory, [80], and as such, they belong to the same kind of general and rather simple construction as their earlier and more particular cases in [1-22,25-32,36-39,41-45,47-49,53,54]. \\

A major interest in the large classes of algebras in the mentioned literature as well as in this paper comes from the following {\it three} properties which are typical and exclusive to these algebras, namely that, their generalized functions

\begin{itemize}

\item are constructed in a rather simple way, requiring only Algebra 101, and specifically, basic results and methods in ring theory, plus a few facts about filters on arbitrary infinite sets,

\item can have singularities on dense subsets of their domain of definition of generalized functions, thus can have sets of singularity points with {\it larger} cardinal than that of their sets of regular, that is, non-singular points, since the only restriction on singularity sets is that their complement, that is, the sets of nonsingular points be dense in the domain of generalized functions,

\item are not subjected to any conditions, in particular, not to growth conditions, in the neighbourhood of their singularities.

\end{itemize}

This nonlinear theory of generalized functions has the further advantage of being given by an {\it infinite} variety of differential algebras of generalized functions, algebras which contain as rather small subspaces the linear space of Schwartz distributions. Consequently, there is a wide liberty in choosing one or another such algebra when dealing with specific problems involving singularities. And not seldom, such a liberty of choice is welcome since it facilitates the appropriate approach to the specific singularities at hand, even if traditionally, one may tend to think that one should rather be given one single algebra, an algebra which would be universally and equally useful in dealing with all possible kind of singularities. \\

As it happens, however, the variety of possible singularities turns out to be so wide as to require more than one single algebra for its proper treatment. Indeed, the essence of this phenomenon is related to the following simple yet fundamental fact

\begin{itemize}

\item the operation of addition does {\it not} appear to branch into alternatives, when dealing with singularities,

\end{itemize}

on the other hand, however

\begin{itemize}

\item the operation of multiplication {\it does} naturally and inevitably branch into infinitely many different alternatives, when dealing with singularities, as shown by most simple algebraic, more precisely, ring theoretic arguments, [7,8,10].

\end{itemize}

The respective algebras do in fact extend, or in other words generalize, large classes of functions $f : X \longrightarrow E$, where $X$ is a domain in an Euclidean space or a finite dimensional manifold, while $E$ is any real or complex, commutative or non-commutative unital Banach algebra, in particular, the usual filed $\mathbb{R}$ of real, or $\mathbb{C}$ of complex numbers. \\

Consequently, these algebras of generalized functions can also deal with non-commutative analysis, in case $E$ are non-commutative Banach algebras. \\

The essence of the applicative utility of these algebras is in the surprisingly large class of singularities which the mentioned nonlinear theory of generalized functions can deal with, and can do so without asking any conditions on generalized functions in the neighbourhood of their singularities. And in an equally surprising manner, much of dealing with such large classes of singularities can be reduced to purely {\it algebraic}, more precisely, ring-theoretic approaches. \\

Here it should be pointed out that the class of admissible singularities is large in no less than {\it two} significantly useful ways : \\

First, the singularities of the functions $f : X \longrightarrow E$ considered can be given by arbitrary subsets $\Sigma \subset X$, subject to the only condition that their complementary $X \setminus \Sigma$, that is, the set of regular, or in other words, non-singular points, be dense in $X$. For instance, if $X = \mathbb{R}^n$ is an Euclidean space, then the set $\Sigma \subset X$ of singularities can be the set of all points with at least one irrational coordinate. Indeed, in this case the set $X \setminus \Sigma$ of non-singular, or regular points is the set of points with all coordinates rational numbers, thus it is dense in $X$. A relevant and rather remarkable fact to note in this case is that the cardinal of the singularity set $\Sigma$ is strictly larger than the cardinal of the set of non-singular points, namely, $X \setminus \Sigma$. \\

Second, there is no restriction on the behaviour of functions $f : X \longrightarrow E$ in the neighbourhood of points in their singularity sets $\Sigma \subset X$, \\

Related to this second freedom in dealing with singularities, one should recall its significant importance in applications. Indeed, as stated in Picard's Great Theorem, an analytic function in the neighbourhood of an isolated singularity point which is an essential singularity takes on all possible complex values infinitely often, with at most a single exception. Consequently, in the neighbourhood of a singularity, one can expect a rather arbitrary behaviour when one deals with more general functions than analytic ones. \\

In this regard, the Colombeau algebras of generalized functions - which are but a particular case of the infinite variety of all possible differential algebras of generalized functions, [8,10-12,15-17,20,22,31,36,40,44] - suffer from a severe limitation. Namely, in the neighbourhood of singularities of their generalized functions, the Colombeau algebras require a polynomial type growth condition, thus they cannot deal even with isolated singularities such as essential singularities of analytic functions. In particular, in Colombeau algebras one simply cannot formulate, let alone obtain, a global version of the classical Cauchy-Kovalavskaia theorem regarding the solution of arbitrary analytic systems of nonlinear PDEs. Indeed, the global solution of such nonlinear PDE systems requires the use of the more powerful differential algebras of generalized functions, algebras which benefit from both above freedoms in dealing with singularities, [4-22,31,32,36,40,44,53,54]. \\
Similarly, due to the essential role of polynomial type growth conditions in the construction of Colombeau algebras, one cannot define arbitrary Lie group actions in Colombeau algebras, actions which require the use of the more powerful algebras in [16,17]. \\

As for the structure of all the differential algebras of generalized functions, they are of the same form of {\it reduced powers}, namely \\

(0.1)~~~ $ A = {\cal A} / {\cal I} $ \\

where \\

(0.2)~~~ $ {\cal I} \subset {\cal A} \subseteq ( {\cal C}^\infty ( X, E ) )^\Lambda $ \\

with $\Lambda$ a suitable infinite set of indices with a directed partial order, while ${\cal A}$ is a subalgebra in $( {\cal C}^\infty ( X, E ) )^\Lambda$, and ${\cal I}$ is an ideal in ${\cal A}$. \\

The fact that the algebras (0.1) are {\it differential} algebras results easily, since the following two conditions can be satisfied in a large variety of situations [6-8,10-17,20,22,31,32,44] \\

(0.3)~~~ $ D^p {\cal A} \subseteq {\cal A},~~~~ p \in \mathbb{N}^n $ \\

(0.4)~~~ $ D^p {\cal I} \subseteq {\cal I},~~~~ p \in \mathbb{N}^n $ \\

in which case the partial derivative operators on the algebras (0.1) can of course be defined by \\

(0.5)~~~ $ A \ni F = f + {\cal I} \longmapsto D^p F = D^p f + {\cal I} \in A,~~~~ p \in \mathbb{N}^n $ \\

where $f = ( f_\lambda ~|~ \lambda\in \Lambda ) \in ( {\cal C}^\infty ( X, E ) )^\Lambda$, while $D^p f = ( D^p f_\lambda ~|~ \lambda\in \Lambda ) \in ( {\cal C}^\infty ( X, E ) )^\Lambda$. \\

The fact that the differential algebras (0.1) - (0.5) contain the Schwartz distributions, thus are algebras of {\it generalized functions} follows from the easy choice of their ideals ${\cal I}$ required to satisfy the {\it off-diagonality} condition, [6-8,10-17,20,22,31,32,44] \\

(0.6)~~~ $ {\cal I} \bigcap {\cal U}_\Lambda = \{ 0 \} $ \\

where ${\cal U}_\Lambda$ is the {\it diagonal} in $( {\cal C}^\infty ( X, E ) )^\Lambda$, that is, the subalgebra of all constant sequences $f = ( f_\lambda ~|~ \lambda\in \Lambda ) \in ( {\cal C}^\infty ( X, E ) )^\Lambda$, where $f_\lambda = f \in {\cal C}^\infty ( X, E )$, for $\lambda \in \Lambda$. \\

And now, in terms of the above, to the motivation of the present paper. \\

As noted, there is a significant interest in {\it enlarging} evermore the classes of singularities that can be dealt with by differential algebras of generalized functions. \\
In this regard, an obvious limitation still present in the general approach in (0.1), (0.2) is the use of ${\cal C}^\infty ( X, E )$, that is, of ${\cal C}^\infty$-smooth functions $f : X \longrightarrow E$. \\

The aim of this paper, therefore, is to replace the use of such ${\cal C}^\infty$-smooth functions in (0.1), (0.2) with the use of a far larger class of functions which are only assumed to be {\it locally} smooth on $X$, see (1.22), (1.23). \\

Needless to say, the classes of singularities which can be dealt with by the resulting differential algebras of generalized functions will not be reduced in any way. \\

As for how large such a class of singularities may be, a rather basic and particular situation can already give a good indication. Namely, in [8,10,11], a {\it global} version was presented for the first time for the classical Cauchy-Kovalevskaia Theorem on the existence of solutions for systems of analytic nonlinear PDEs. And the respective global existence result shows the extent to which the generalized solutions must typically be {\it local}, even if where they exist, they turn out to be {\it analytic}. \\
Clearly, therefore, in the general, non-analytic case of nonlinear systems of PDEs, the typical generalized solutions cannot be expected to have simpler singularities. \\
As for how complex the structure of such singularities may be, relevant recent information can be found in [33-35,40,62-64], an information which, however, is at present not yet explicit enough. \\

Consequently, the interest in the differential algebras of generalized functions constructed in this paper, algebras which allow for significantly more general singularities than the earlier similar ones constructed in [1-22,31,32,36,44,54]. \\ \\

{\bf 1. Basic Definitions} \\

{\bf A Single Singularity Set} \\

Let $X$ be a nonvoid set and $\Sigma \subset X$ which will play the role of the subset of {\it singular} points of certain functions defined locally on $X$. \\

{\bf Definition 1.1.} \\

A $\Sigma$-{\it local function} on $X$ is every family \\

(1.1)~~~ $ f = ( f_{x, \, U_x} ~|~ x \in X \setminus \Sigma ) $ \\

where $x \in U_x \subseteq X$, while $f_{x, \, U_x} : U_x \longrightarrow E$, and the following compatibility condition holds \\

(1.2)~~~ $ \forall~ x, y \in X \setminus \Sigma ~ : ~ x \in U_y, \, y \in U_x ~~\Longrightarrow~~
                                                  f_{x, \, U_x} = f_{y, \, U_y} ~\mbox{on}~ U_x \cap U_y $ \\

{\bf Remark 1.1.} \\

We note that condition (1.2) is considerably weaker than condition \\

(1.3)~~~ $ \forall~ x, y \in X \setminus \Sigma ~ : ~ U_x \cap U_y \neq \phi ~~\Longrightarrow~~
                                                  f_{x, \, U_x} = f_{y, \, U_y} ~\mbox{on}~ U_x \cap U_y $ \\

as seen in 2) in the following \\

{\bf Examples 1.1.} \\

1) Let $f : X \longrightarrow \mathbb{R}$, then $( f|_{ \, U_x} ~|~ x \in X \setminus \Sigma )$ is a $\Sigma$-local function on $X$, whenever $U_x \subseteq X$ are a neighbourhoods of $x \in X$. \\

2) Let $X = \mathbb{R}, \, \Sigma = \mathbb{R} \setminus \mathbb{Q}$, where as usual, $\mathbb{Q}$ denotes the set of rational numbers. We assume $X \setminus \Sigma = \{ x_0, x_1, x_2, \ldots \}$ and take \\

$U_0 = ( x_0 - r_0, x_0 + r_0)$ with $r_0 > 0$ \\

$U_1 = ( x_1 - r_1, x_1 + r_1)$ with $r_1 > 0$, such that $x_0 \notin U_1$ \\

$U_2 = ( x_2 - r_2, x_2 + r_2)$ with $r_2 > 0$, such that $x_0, x_1 \notin U_2$ \\

$U_3 = ( x_3 - r_3, x_3 + r_3)$ with $r_3 > 0$, such that $x_0, x_1, x_2 \notin U_3$ \\

$\vdots$ \\

Given now $c_0, c_1, c_2, \ldots \in \mathbb{R}$, we define $f_{x_n, \, U_{x_n}} = c_n$, for $n \in \mathbb{N}$. Then \\

(1.4)~~~ $ f = ( f_{x, \, U_x} ~|~ x \in X \setminus \Sigma )$ is a  $\Sigma$-local function on $X$ \\

Indeed, let $n < m$. Then $x_0, x_1, \ldots , x_{m - 1} \notin U_m$, hence $x_n \notin U_m$, therefore (1.2) is satisfied by default. \\

Clearly, for $n \in \mathbb{N}$, there are infinitely many $m \in \mathbb{N}$, such that \\

(1.5)~~~ $ U_n \cap U_m \neq \phi $ \\

therefore, condition (1.3) is in general not satisfied. \\

The interest in (1.4) is in the following five facts : first, the set $X \setminus \Sigma$ is dense in $X$, second, the values $c_0, c_1, c_2, \ldots \in \mathbb{R}$ can be arbitrary, third, the sum $\sum_{n \in \mathbb{N}} r_n$ can be arbitrary small, thus so can be the measure of $\bigcup_{n \in \mathbb{N}} U_n$, fourth, the sets $U_x$ are open neighbourhoods of the respective $x$, fifth, the component functions $f_{x_n, \, U_{x_n}}$ are highly smooth or regular, being in fact constant.

\hfill $\Box$ \\

We denote by \\

(1.6)~~~ $ {\cal B}_{{lc}, \, \Sigma} ( X, E ) $ \\

the set of all $\Sigma$-local functions on $X$ with values in the Banach algebra $E$. Clearly, ${\cal B}_{{lc}, \, \Sigma} ( X, E )$  is a commutative, respectively, non-commutative unital algebra on $\mathbb{R}$, according to $E$ being commutative or not. Further, we define the {\it algebra embedding} \\

(1.7)~~~ $ E^X \ni f \longrightarrow
                  lc ( f ) = ( f_{x, \, X} ~|~ x \in X \setminus \Sigma ) \in {\cal B}_{{lc}, \, \Sigma} ( X, E ) $ \\

where $f_{x, \, X} = f$. Also, we denote by \\

(1.8)~~~ $ {\cal V}_{{lc}, \, \Sigma} ( X, E ) $ \\

the subalgebra in ${\cal B}_{{lc}, \, \Sigma} ( X, E )$ which is the {\it range} of the above algebra embedding (1.7).
Thus we have the {\it algebra isomorphism} \\

(1.9)~~~ $ E^X \ni f \longrightarrow
                  lc ( f ) = ( f_{x, \, X} ~|~ x \in X \setminus \Sigma ) \in {\cal V}_{{lc}, \, \Sigma} ( X, E ) $ \\

Given $Z \subseteq X$, with $Z \setminus \Sigma \neq \phi$, we denote by \\

(1.10)~~~ $ {\cal J}_{{lc}, \, \Sigma, \ Z} ( X, E ) $ \\

the set of all $f = ( f_{ x, \, U_x } ~|~ x \in X \setminus \Sigma ) \in {\cal B}_{{lc}, \, \Sigma} ( X, E )$, such that \\

(1.11)~~~ $ \forall~~ x \in Z \setminus \Sigma ~:~ f_{ x, \, U_x } ( x ) = 0 $ \\

Obviously, ${\cal J}_{{lc}, \, \Sigma, \ Z} ( X, E )$ is an {\it ideal} in ${\cal B}_{{lc}, \, \Sigma} ( X, E )$. \\

Given now $\Sigma \subseteq \Sigma' \subset X$, we can define the mapping \\

(1.12)~~~ $ j_{\, \Sigma, \, \Sigma'} : {\cal B}_{{lc}, \, \Sigma} ( X, E )
                                  \longrightarrow {\cal B}_{{lc}, \, \Sigma\,'} ( X, E )$ \\

by \\

(1.13)~~~~ $ j_{\, \Sigma, \, \Sigma'} ( f_{x, \, U_x} ~|~ x \in X \setminus \Sigma ) =
                                           ( f_{x, \, U_x} ~|~ x \in X \setminus \Sigma' ) $ \\

And these mappings are {\it surjective algebra homomorphisms} which have the properties \\

(1.14)~~~ $ j_{\, \Sigma, \, \Sigma} = id_{{\cal B}_{{lc}, \, \Sigma} ( X, E )} $, for $ \Sigma \subset X $ \\

(1.15)~~~ $j_{\, \Sigma', \, \Sigma''} \circ j_{\, \Sigma, \, \Sigma'} = j_{\, \Sigma, \, \Sigma'' } $, for
                           $ \Sigma \subseteq \Sigma' \subseteq \Sigma'' \subset X $ \\

{\bf Families of Singularity Sets} \\

Let us now turn to the case when instead of one single subset $\Sigma \subset X$ of singularities, we have a whole family ${\cal S} \subseteq {\cal P} ( X )$ of such singularity subsets $\Sigma \in {\cal S}$. In this regard, we shall assume in the sequel that \\

(1.16)~~~ $ X \notin {\cal S} $ \\

(1.17)~~~ $ \forall~~ \Sigma, \Sigma' \in {\cal S} ~:~ \exists~~ \Sigma'' \in {\cal S} ~:~
                                                           \Sigma \cup \Sigma' \subseteq \Sigma'' $ \\

Obviously, (1.16) is equivalent with $\Sigma \subset X$, for $\Sigma \in {\cal S}$. \\

In view of (1.17), it follows that $( {\cal S}, \subseteq )$ is a {\it directed} partially ordered set. \\

Clearly, in the particular case when ${\cal S} = \{ \Sigma \}$, that is, when we have one single subset $\Sigma \subset X$ of singularities, then the conditions (1.16), (1.17) are satisfied. \\

We consider now in the general case of (1.16), (1.17), the set \\

(1.18)~~~ $ {\cal B}_{{lc}, \, {\cal S}} ( X, E ) =
                    \bigcup_{\, \Sigma \in \, {\cal S}} \, {\cal B}_{{lc}, \, \Sigma} ( X, E ) $ \\

as well as \\

(1.19)~~~~ $ {\cal V}_{{lc}, \, {\cal S}} ( X, E ) =
                                \bigcup_{\, \Sigma \in \, {\cal S}} {\cal V}_{{lc}, \, \Sigma} ( X, E ) $ \\

which is obviously a subset of ${\cal B}_{{lc}, \, {\cal S}} ( X, E )$. \\

Further, let $Z \subseteq X$, such that \\

(1.20)~~~ $ \forall~~ \Sigma \in {\cal S} ~:~ Z \setminus \Sigma \neq \phi$ \\

Then we define \\

(1.21)~~~ $  {\cal J}_{{lc}, \, {\cal S}, \, Z} ( X, E ) =
                    \bigcup_{\, \Sigma \in \, {\cal S}} \, {\cal J}_{{lc}, \, \Sigma, \, Z} ( X, E )$ \\

which is obviously a subset of ${\cal B}_{{lc}, \, {\cal S}} ( X, E )$. \\

Now we consider the {\it direct limit} \\

(1.22)~~~ $ {\cal A}_{\, {lc}, \, {\cal S}} ( X, E ) =
                     \underrightarrow{\lim}_{\, \Sigma \in \, {\cal S}} \, {\cal B}_{{lc}, \, \Sigma} ( X, E ) $ \\

It follows that \\

(1.23)~~~ $ {\cal A}_{\, {lc}, \, {\cal S}} ( X, E ) = {\cal B}_{{lc}, \, {\cal S}} ( X, E ) / \approx_{\cal S} $ \\

where the equivalence relation $\approx_{\cal S}$ on ${\cal B}_{{lc}, \, {\cal S}} ( X, E )$ is defined for $ ( f_{x, \, U_x} ~|~ x \in X \setminus \Sigma ) \in {\cal B}_{{lc}, \, \Sigma} ( X, E ), \, ( g_{y, \, V_y} ~|~ y \in X \setminus \Sigma' ) \in {\cal B}_{{lc}, \, \Sigma\,'} ( X, E )$, with $\Sigma, \Sigma' \in {\cal S}$, by \\

(1.24)~~~ $ ( f_{x, \, U_x} ~|~ x \in X \setminus \Sigma ) \approx_{\cal S}
                                     ( g_{y, \, V_y} ~|~ y \in X \setminus \Sigma' ) $ \\

if and only if there exist $\Sigma'' \in {\cal S}$, with $\Sigma \cup \Sigma' \subseteq \Sigma''$, as well as $( h_{z, \, W_z} ~|~ z \in X \setminus \Sigma'' ) \in {\cal B}_{{lc}, \, \Sigma''} ( X, E ) $, such that $j_{\, \Sigma, \, \Sigma''} ( f_{x, \, U_x} ~|~ x \in X \setminus \Sigma ) = j_{\, \Sigma', \, \Sigma''} ( g_{y, \, V_y} ~|~ y \in X \setminus \Sigma' ) =
( h_{z, \, W_z} ~|~ z \in X \setminus \Sigma'' )$. \\

Similarly, one defines the {\it direct limit} \\

(1.25)~~~ $ {\cal U}_{\, lc, \, {\cal S}} ( X, E ) =
                \underrightarrow{\lim}_{\, \Sigma \in \, {\cal S}} \, {\cal V}_{{lc}, \Sigma} ( X, E ) $ \\

and obtains \\

(1.26)~~~ $ {\cal U}_{\, lc, \, {\cal S}} ( X, E ) =
                    {\cal V}_{\, lc, \, {\cal S}} ( X, E ) ) / \approx_{\cal S} $ \\

Clearly, we have the {\it injective} mapping \\

(1.27)~~~ $ E^X \ni f \longrightarrow
                  ( lc ( f ) )_ {\approx_{\cal S}} \in {\cal A}_{\, lc, \, {\cal S}} ( X, E ) $ \\

where $( g )_ {\approx_{\cal S}}$ denotes the $\approx_{\cal S}$ equivalence class of the element $g \in {\cal B}_{\, lc, \, {\cal S}} ( X, E )$. \\

We note that, if ${\cal S} = \{ \Sigma \}$, then ${\cal A}_{\, lc, \, {\cal S}} ( X, E ) = {\cal A}_{\, {lc}, \, \Sigma} ( X, E ) = {\cal B}_{\, {lc}, \, \Sigma} ( X, E )$ and ${\cal U}_{\, lc, \, {\cal S}} ( X, E ) = {\cal U}_{\, lc, \, \Sigma} ( X, E ) = {\cal V}_{\, lc, \, \Sigma} ( X, E )$. \\

Lastly, given $Z \subseteq X$ for which (1.20) holds, one defines the {\it direct limit} \\

(1.28)~~~ $ {\cal I}_{\, lc, \, {\cal S}, \, Z} ( X, E ) =
               \underrightarrow{\lim}_{\, \Sigma \in \, {\cal S}} \, {\cal J}_{{lc}, \, \Sigma, \, Z} ( X, E ) $ \\

and obtains \\

(1.29)~~~ $ {\cal I}_{\, lc, \, {\cal S}, \, Z} ( X, E ) =
                    {\cal J}_{\, lc, \, {\cal S}, \, Z} ( X, E ) ) / \approx_{\cal S} $ \\

Obviously, in view of (1.21), we have \\

(1.30)~~~ $  {\cal I}_{\, lc, \, {\cal S}, \, Z} ( X, E ) \subseteq  {\cal A}_{\, lc, \, {\cal S}} ( X, E ) $ \\

We can also note that, if ${\cal S} = \{ \Sigma \}$, then ${\cal I}_{\, lc, \, {\cal S}, \, Z} ( X, E ) = {\cal I}_{\, lc, \, \Sigma, \, Z} ( X, E ) = {\cal J}_{\, lc, \, \Sigma, \, Z} ( X, E )$ \\

Let us summarize. Given $\Sigma \in {\cal S}$ as above, and $Z \subseteq X$ as in (1.20), we have the {\it commutative} diagram of mappings \\

$~~ \begin{array}{l}
        E^X \stackrel{*}\longrightarrow {\cal V}_{lc, \, \Sigma} ( X, E ) \longrightarrow {\cal B}_{lc, \, \Sigma}
             ( X, E ) \longleftarrow {\cal J}_{\, lc, \, \Sigma, \, Z} ( X, E ) \\
        \hspace*{2cm} \downarrow \hspace*{2.7cm} \downarrow \hspace*{3.2cm} \downarrow  \\
        \hspace*{1.6cm} {\cal V}_{lc, \, {\cal S}} ( X, E ) \longrightarrow
             {\cal B}_{lc, \, {\cal S}} ( X, E ) \longleftarrow {\cal J}_{\, lc, \, {\cal S}, \, Z} ( X, E ) \\
        \hspace*{2cm} \downarrow* \hspace*{2.3cm} \downarrow* \hspace*{2.8cm} \downarrow* \\
        E^X \stackrel{*}\longrightarrow {\cal U}_{lc, \, {\cal S}} ( X, E ) \longrightarrow {\cal A}_{lc, \, {\cal S}}
             ( X, E ) \longleftarrow {\cal I}_{\, lc, \, {\cal S}, \, Z} ( X, E ) \\
    \end{array} $ \\ \\

where all the mappings are {\it injective}, except for the three mappings ''$\downarrow*$'' which are {\it surjective}, while the two mappings ''$\stackrel{*}\longrightarrow$'' are in fact {\it bijective}. \\ \\

{\bf 2. Properties} \\

The following result can be obtained by direct, even if a somewhat elaborate verification : \\

{\bf Theorem 2.1.} \\

${\cal A}_{\, {lc}, \, {\cal S}} ( X, E )$ is a unital algebra, and ${\cal U}_{\, lc, {\cal S}} ( X, E )$ is a subalgebra in it, and it is the range of the mapping (1.27) which is an {\it algebra embedding}, namely, we have the {\it algebra isomorphism}  \\

(2.1)~~~ $ E^X \ni f \longrightarrow
                  ( lc ( f ) )_ {\approx_{\cal S}} \in {\cal U}_{\, lc, \, {\cal S}} ( X, E ) \subset {\cal A}_{\, lc, \, {\cal S}} ( X, E ) $ \\

The algebra ${\cal A}_{\, {lc}, \, {\cal S}} ( X, E )$ is commutative, if and only if the Banach algebra $E$ is commutative. \\

Furthermore, given $Z \subseteq X$ for which (1.20) holds, then ${\cal I}_{\, lc, \, {\cal S}, \, Z} ( X, E )$ is an {\it ideal} in the algebra ${\cal A}_{\, {lc}, \, {\cal S}} ( X, E )$. \\

\hfill $\Box$ \\

For $0 \leq l \leq \infty$, we denote by, see (1.22), (1.23) \\

(2.2)~~~ $ {\cal A}^l_{\, {lc}, \, {\cal S}} ( X, E ) $ \\

the set of all $( f )_{\approx_{\cal S}}$ where $f = ( f_{x, \, U_x} ~|~ x \in X \setminus \Sigma ) \in {\cal B}_{\, lc, \, \Sigma} ( X, E )$, for some $\Sigma \in {\cal S}$, such that \\

(2.3)~~~ $ f_{x, \, U_x} \in {\cal C}^l ( U_x ) $, \, for \, $ x \in X \setminus \Sigma $ \\

Further, we denote \\

(2.4)~~~ $  {\cal U}^{\, l}_{\, lc, \, {\cal S}} ( X, E ) = {\cal U}_{\, lc, \, {\cal S}} ( X, E ) \, \bigcap \,
                                   {\cal A}^l_{\, {lc}, \, {\cal S}} ( X, E ) $ \\

(2.5)~~~ $  {\cal I}^{\, l}_{\, lc, \, {\cal S}, \, Z} ( X, E ) = {\cal I}_{\, lc, \, {\cal S}, \, Z} ( X, E ) \, \bigcap \,
                                   {\cal A}^l_{\, {lc}, \, {\cal S}} ( X, E ) $ \\

{\bf Theorem 2.2.} \\

Suppose that \\

(2.6)~~~ $ \forall~~ \Sigma \in {\cal S} ~:~ X \setminus \Sigma $ \, is dense in \, $ X $ \\

 Then, for $0 \leq l \leq \infty$, the ideal ${\cal I}^{\, l}_{\, lc, \, {\cal S}} ( X, E )$ in the algebra
 ${\cal A}^l_{\, {lc}, \, {\cal S}} ( X, E )$ satisfies the {\it off-diagonality} condition \\

(2.7)~~~ $ {\cal I}_{\, lc, \, {\cal S}} ( X, E ) \, \bigcap \,
                           {\cal U}_{\, lc, \, {\cal S}} ( X, E ) = \{ 0 \} $ \\

{\bf Proof} \\

It follows from the density condition (2.6) and the continuity of the functions involved. Indeed, let, see (1.29) \\

(2.8)~~~ $ ( f )_{\approx_{\cal S}} \in {\cal I}_{\, lc, \, {\cal S}} ( X, E ) $ \\

 where for some $\Sigma \in {\cal S}$, we have $f = ( f_{ x, \, U_x } ~|~ x \in X \setminus \Sigma ) \in {\cal B}_{lc, \, \Sigma} ( X, E )$. Then $( f )_{\approx_{\cal S}} \in {\cal U}_{\, lc, \, {\cal S}} ( X, E )$ implies in view of (1.26) that, see (1.8) \\

(2.9)~~~ $ f = ( f_{ x, \, U_x } ~|~ x \in X \setminus \Sigma ) \in {\cal V}_{{lc}, \, \Sigma} ( X, E ) $ \\

Now (2.9), (1.8) give \\

(2.10)~~~ $ f_{ x, \, U_x } = f,~~~ x \in X $ \\

for some $f \in {\cal C}^0 ( X, E )$. \\

On the other hand, (2.8), (1.11), (2.10) give \\

$~~~~~~ f ( x ) = f_{ x, \, U_x } ( x ) = 0,~~~ x \in X $ \\

thus indeed $f = 0$. \\ \\

{\bf 3. Differential Algebras with Dense Singularities} \\

In [18-21,31,44], large classes of differential algebras of generalized functions which allow their elements to have singularities on {\it dense} subsets of their domain of definition, and {\it without} any restrictions on the respective generalized functions in the neighbourhood of their singularities, have been introduced, and these algebras have been applied to solving large classes of systems of nonlinear PDEs, as well as in highly singular problems in Differential Geometry, with respective applications in modern Physics, including General Relativity and Quantum Gravity, [9,13-19,21,56-61]. \\

The algebras are of the form (0.1), (0.2), thus are built upon ${\cal C}^\infty$-smooth functions, namely, their elements are classes of equivalence modulo the respective ideals ${\cal I}$. \\

These ideals play a fundamental role in dealing with dense singularities, and so without any restrictions on the respective generalized functions in the neighbourhood of their singularities. Indeed, the power of the method consists precisely in the fact that the singularities, although possibly so many as to constitute dense subsets in the domain of definition of generalized functions, are dealt with exclusively {\it algebraic}, that is, ring theoretic means. \\
And here it should be mentioned that the singularities can forms sets which have a {\it larger} cardinal then the set of regular, that is, non-singular points. For instance, if the generalized functions are defined on $X = \mathbb{R}$, then the set of singularities can be given by all irrational numbers, thus the set of regular, non-singular points can be reduced to the set of rational numbers, [18-21,31,44]. \\

As argued in section 0, there is a major interest in extending such algebras by replacing the ${\cal C}^\infty$-smooth functions upon which they are built with considerably larger classes of functions, and specifically in this paper, with functions which are {\it locally} smooth, see (1.22), (1.23). \\

When proceeding with such an extension, the main issue is to extend in appropriate ways the definition of the large class of ideals ${\cal I}$ in such a way that they still can handle {\it dense} singularities, and do so {\it without} any restrictions on the respective generalized functions in the neighbourhood of their singularities, just as they were able to do in [18-21,31,44]. \\

Let us, therefore, recall for convenience the definition of the large class of ideals ${\cal I}$ in [18-21,31,44], in the case when, as assumed, $X$ is a domain in and Euclidean space $\mathbb{R}^n$. \\

First we recall that we assumed a directed partial order $\leq$ on the infinite set of indices $\Lambda$. Further, let, as in (1.16), (1.17), ${\cal S} \subseteq {\cal P} ( X )$ be a family of singularity subsets $\Sigma \subset X$. \\

Now for given $\Sigma \in {\cal S}$, we considered in [18-21,31,44] the ideal \\

(3.1)~~~ $ {\cal I}_\Sigma ( X ) $ \\

in $( {\cal C}^\infty ( X ) )^\Lambda$, given by all the sequences of smooth functions $w = ( w_\lambda ~|~ \lambda \in \Lambda ) \in ( {\cal C}^\infty ( X ) )^\Lambda$, such that \\

(3.2)~~~ $ \begin{array}{l}
               \forall~~ x \in X \setminus \Sigma ~: \\
               \exists~~ \lambda \in \Lambda ~: \\
               \forall~~ \mu \in \Lambda, \, \mu \geq \lambda ~: \\
               \forall~~ p \in \mathbb{N}^n ~: \\
               ~~~~ D^p w_\mu ( x ) = 0
            \end{array} $ \\

Further, we defined the ideal in $( {\cal C}^\infty ( X ) )^\Lambda$, given by \\

(3.3)~~~ $ {\cal I}_{\cal S} ( X ) = \bigcup_{\, \Sigma \in \, {\cal S}} {\cal I}_\Sigma ( X ) $ \\

which played the role of the ideals ${\cal I}$ in (0.1), (0.2). \\

In order to extend these ideals to the case of {\it locally} smooth functions, first we extend (3.2) according to (1.1) - (1.7). Namely, we denote by \\

(3.4)~~~ $ {\cal N}_{lc, \, {\cal S}} ( X, E ) $ \\

the set of all the sequences $\widetilde w$, where for a suitable $\Sigma \in {\cal S}$, we have $\widetilde w = ( ( w_\lambda )_{\approx_{\cal S}} ~|~ \lambda \in \Lambda ) \in ( {\cal B}^\infty_{lc, \, \Sigma} ( X, E ) )^\Lambda$, such that \\

(3.5)~~~ $ \begin{array}{l}
               \forall~~ x \in X \setminus \Sigma ~: \\
               \exists~~ \lambda \in \Lambda ~: \\
               \forall~~ \mu \in \Lambda, \, \mu \geq \lambda ~: \\
               \forall~~ p \in \mathbb{N}^n ~: \\
               ~~~~ D^p w_\mu ( x ) = 0
            \end{array} $ \\ \\

{\bf Theorem 3.1.} \\

${\cal N}_{lc, \, {\cal S}} ( X, E )$ is an ideal in $( {\cal A}^\infty_{lc, \, {\cal S}} ( X, E ) )^\Lambda$. \\

{\bf Proof.} \\

It follows by direct verification, based on (2.2), (2.3) and (3.4), (3.5).

\hfill $\Box$ \\

At last, we can now arrive at the {\it main construction} in this paper, namely, the {\it reduced power algebras} \\

(3.6)~~~ $ A_{lc, \, {\cal S}} ( X, E ) =
                  ( {\cal A}^\infty_{lc, \, {\cal S}} ( X, E ) )^\Lambda / {\cal N}_{lc, \, {\cal S}} ( X, E ) $ \\

which are in fact {\it differential algebras of generalized functions}. Furthermore, they contain all the earlier differential algebras of generalized functions, [1-17,22], and in particular, those with dense singularities, [18-21], the Colombeau algebras, and therefore, also the linear vector spaces of Schwartz distributions, [6-8,10]. \\

Indeed, it follows from a direct, albeit elaborate verification that the algebras (3.6) satisfy the corresponding conditions (0.3) - (0.6). \\ \\

{\bf 4. Comments} \\

1) The Model Theoretic, [80], construction of {\it reduced power}, although hardly known as such among so called working mathematicians, happens nevertheless to appear in quite a number of important places in Mathematics at large. For a sample of them, one can note the following. The Cauchy-Bolzano construction of the field $\mathbb{R}$ of usual real numbers is in fact a reduced power of the rational numbers $\mathbb{Q}$. More generally, the completion of any metric space is a reduced power of that space. Furthermore, this is but a particular case of the fact that the completion of any uniform topological space is a reduced power of that space. Also, in a rather different direction, the field $^*\mathbb{R}$ of nonstandard real numbers can be obtained as a reduced power of the usual field $\mathbb{R}$ of real numbers. \\

In view of the above, the use of reduced powers in the construction of differential algebras of generalized functions should not be seen as much else but a further application of that basic construction in Model Theory, this time to the study of large classes of singularities. \\

As for dealing with singularities, there is a strongly entrenched trend to approach them with nothing else but methods of Analysis, Functional Analysis, Topology, or Complex Functions. And this trend is particularly manifest in various theories of generalized functions. \\

On the other hand, as seen in [1-22,25-32,36-39,41-45,47-49,53,54], and specifically, in [10], the issue of singularities of generalized functions boils down to a rather simple and basic {\it algebraic conflict}. Consequently,
the study of singularities of generalized functions through methods which are primarily of Analysis, Functional Analysis, Topology, or Complex Functions has the {\it double disadvantage} of

\begin{itemize}

\item unnecessarily complicating the situation,

\item missing the root of the problem.

\end{itemize}

The Colombeau algebras do to a good extent fall under the above double disadvantage, and as result, mentioned in section 0, they can only deal with a rather limited class of singularities of generalized functions. \\

2) An important property of many reduced powers is the presence of {\it infinitesimals}, [52]. This fact, as it happens, has not yet been given its due consideration in the study of differential algebras of generalized functions. \\

3) A large class of {\it scalars} which most likely may have a considerable relevance in Physics is given by reduced power algebras built upon the field $\mathbb{R}$ of usual real numbers. Indications in this regard can be found in [25,26,37-39,41-43,45,47-49]. \\

\end{document}